\documentclass[a4paper,11pt,bibliography=totoc,twoside,reqno]{amsart}

\usepackage[utf8]{inputenc}

\usepackage{amssymb,amsthm}
\usepackage[plainpages=false,pdfpagelabels=true]{hyperref}
\usepackage[margin=1in]{geometry}

\newtheorem{Satz}{Theorem}[section]
\newtheorem{Prop}[Satz]{Proposition}
\newtheorem{Lem}[Satz]{Lemma}
\newtheorem{Thm}[Satz]{Theorem}
\newtheorem{Cor}[Satz]{Corollary}
\newcommand{\vol}{{\operatorname{Vol}}}
\theoremstyle{definition}

\newtheorem{Bem}[Satz]{Remark}

\allowdisplaybreaks[1]

\renewcommand{\epsilon}{\varepsilon}

\newcommand{\R}{\ensuremath{\mathbb{R}}}

\newcommand{\Z}{\ensuremath{\mathbb{Z}}}
\newcommand{\N}{\ensuremath{\mathbb{N}}}

\newcommand{\hyp}{\ensuremath{\mathbb{H}}}
\newcommand{\cG}{{\mathcal G}}

\newcommand{\hol}{\operatorname{Hol}}
\newcommand{\Hom}{\operatorname{Hom}}

\numberwithin{equation}{section}

\title{Magnetic Geodesics via the Heat Flow}
\author{Volker Branding and Florian Hanisch}
\date{\today}
\address{TU Wien\\
Institut für diskrete Mathematik und Geometrie\\
Wiedner Hauptstraße 8–10, A-1040 Wien}
\email[]{volker@geometrie.tuwien.ac.at}
\address{Universität Potsdam\\
Institut für Mathematik\\
Am Neuen Palais 10, 14469 Potsdam}
\email[]{fhanisch@math.uni-potsdam.de}
\begin{document}

\begin{abstract}
Magnetic geodesics describe the trajectory of a particle in a Riemannian manifold under the influence of an external magnetic field.
In this article, we use the heat flow method to derive existence results for such curves. 
We first establish subconvergence of this flow to a magnetic geodesic under certain boundedness assumptions. 
It is then shown that these conditions are satisfied provided that either the magnetic field admits a global potential or the initial curve is sufficiently small. 
Finally, we discuss different examples which illustrate our results.
\end{abstract} 

\maketitle

\section{Introduction and Results}
The trajectory of a free point particle in a Riemannian manifold is modeled by a geodesic.
Geodesics arise from a variational problem, the existence of closed solutions can be studied by various methods and is well understood at present. 
A powerful method to ensure their existence is the so called heat flow method. Here, one deforms a given initial curve \(\gamma_0\) by a heat-type equation and obtains a closed geodesic in the end.
For geodesics, this method was successfully applied by Ottarsson \cite{MR834094} building on
the famous existence result for harmonic maps into manifolds with non-positive curvature due to Eells and Sampson \cite{MR0164306}.
\par\medskip
The aim of this article is to study the equation which determines the trajectory of a point particle in a Riemannian manifold
in the presence of an external magnetic field. This equation is known as equation for \emph{magnetic geodesics} or \emph{prescribed geodesic curvature equation}.
This equation can also be derived from a variational principle, however the corresponding energy functional is \(U(1)\)-valued in general.
\par\medskip
The existence of magnetic geodesics has been ensured by various methods.
This includes methods from dynamical systems and symplectic geometry \cite{MR890489}, \cite{MR902290}, \cite{MR1458315}, \cite{MR1888853}, \cite{MR1432462}, \cite{MR1417851}, \cite{MR2250797}
and variational methods for multivalued functionals (Morse-Novikov-theory) \cite{MR1133303}, \cite{MR1185286}, \cite{MR730159}.
Finally, one can also derive an existence result by Aubry-Mather's theory \cite{MR2036336}.
Recently, a new approach was introduced: An existence result for magnetic geodesics on \(S^2\) \cite{MR2788659}
and also on closed hyperbolic surfaces \cite{MR2959932} is established by studying the zeros of a certain vector field.
\par\medskip
In this article we use the heat flow method to approach the existence of magnetic geodesics.
Namely, we deform a given initial curve \(\gamma_0\) by a heat type equation
and study in which cases this equation converges to magnetic geodesic. This approach was initiated
in \cite{MR2551140}.
In \cite{MR1800592} a similar problem is discussed, namely the heat flow
for harmonic maps coupled to a potential.

Let us describe the problem in more detail. Suppose that $N$ is a closed Riemannian manifold and \(\gamma\colon S^1\to N\) is a smooth curve.
The $U(1)$-valued  energy functional for magnetic geodesics is given by
\begin{equation}
\label{energy-functional-holonomy}
C^\infty(S^1,N)\to U(1),\qquad\gamma\mapsto e^{i\frac{1}{2}\int_{S^1}|\gamma'|^2ds}\hol(\gamma), 
\end{equation}
where \(\hol(\gamma)\) represents the holonomy of the magnetic field (more precisely, of the corresponding gerbe \(\cG\)) along \(\gamma\). 
Moreover, \(\gamma'\) represents the derivative with respect to the curve parameter \(s\).

In the case that the magnetic field is given by an exact form (or in the language of gerbes: the curvature 2-form is exact),
we may rewrite the energy functional as the sum of a kinetic and a magnetic contribution 
\begin{equation}
\label{energy-functional}
E(\gamma)= E_{kin}(\gamma) + E_A(\gamma) := \frac{1}{2}\int_{S^1}|\gamma'|^2ds+\int_{S^1}\gamma^\ast Ads,
\end{equation}
where the one-form \(A\) is the potential for the magnetic field.
The critical points of \eqref{energy-functional-holonomy} and \eqref{energy-functional} are given by
\begin{equation}
\tau(\gamma)=Z(\gamma'),
\end{equation}
where \(\tau\) denotes the tension field of the curve \(\gamma\) and \(Z\in\Gamma (\Hom(TN,TN))\) represents the magnetic field strength. 
Note that the critical points are globally defined, even in the case that the energy functional is not.

In the following, we are studying the \(L^2\)-gradient flow of the energy functional (\ref{energy-functional}),
which is given by
\begin{equation}
\label{gradient-flow}
\left\{
  \begin{array}{l l}
    \dot{\gamma}_t(s,t)=\tau(\gamma_t)(s,t)-Z(\gamma'_t)(s,t),\qquad (s,t)\in S^1\times [0,\infty),\\
    \gamma(s,0)=\gamma_0(s)\\
  \end{array} \right.
\end{equation}
and we use a \(\dot{\gamma}\) to denote the derivative of $\gamma$ with respect to \(t\). Note that neither the functionals \eqref{energy-functional-holonomy} and \eqref{energy-functional} 
nor the  equations of motion derived from them are invariant under rescaling of the domain of definition. Thus, strictly speaking, we have a different problem for each class of curves $C^\infty(S^1_r, N)$, 
where $S^1_r$ denotes the circle of length $2\pi r$. 
From a physical point of view, the length of the circle encodes the period of revolution of the charged particle. 
Since we are interested in the existence of closed trajectories regardless of this time, we ultimately look for critical points defined on some $S^1_r$. 
We will nevertheless mostly work with $S^1 := S^1_1$ but make use of the possibility of rescaling the circle in Corollary \ref{CorRescaling}.\\

The existence of a short-time solution of the evolution equation \eqref{gradient-flow} with existence interval \([0,T)\) 
for initial data \(\gamma_0\in C^{2+\alpha}(S^1,N)\) is guaranteed by Theorem 19 in \cite{MR2551140}.
Moreover, Theorem 1 in the same reference shows that \eqref{gradient-flow} has a unique, smooth solution for all \(t\in[0,\infty)\). 
However, the question in which cases the gradient flow converges was not fully answered in \cite{MR2551140}.\\

In this article, we will prove the following main results:
\begin{Satz}
Let \(\gamma_t\colon S^1\times [0,\infty)\to N\) denote the unique solution of \eqref{gradient-flow} associated to \(Z\in\Gamma (\Hom(TN,TN))\) and initial condition \(\gamma_0\in C^{2+\alpha}(S^1,N)\). 
\begin{enumerate}
 \item If the magnetic field admits a global potential, then (\ref{gradient-flow}) subconverges and we obtain a smooth magnetic geodesic \(\gamma_\infty\).
       If $\gamma_0$ is not null-homotopic or $E(\gamma_0) \leq 0$ and $\gamma_0$ is not constant, then the magnetic geodesic \(\gamma_\infty\) is not trivial (i.e. not a point).
 \item If $|Z|_{L^\infty}$ is sufficiently small and the initial curve has sufficiently small kinetic energy (see Lemma \ref{LemApplyOttarson} for a precise formulation), then (\ref{gradient-flow}) 
       subconverges. Under a slightly stronger assumption on $|Z|_{L^\infty}$, $\gamma_t$ converges to a point. \label{Item2}
\end{enumerate}
\end{Satz}
Case \eqref{Item2}, where the magnetic field is not derived from a global potential,
is particularly challenging since \eqref{gradient-flow} is no longer the gradient flow associated to a globally defined energy.\\

This paper is organized as follows: In Section 2 we recall the derivation of the equation for magnetic geodesics
and the Bochner formulas for the associated evolution equation. Section 3 then discusses the convergence of
the gradient flow. We give a general criterion for subconvergence of the flow (Theorem \ref{Convergence}) 
and show that it can be applied to the cases where either the field admits a global potential (Theorem \ref{ThmConvExactCase})
or the initial curve has sufficiently low energy and the magnetic field is weak (Theorem \ref{ThmConvergenceOttarson}). 
Finally, in Section 4 we calculate some examples and compare our general theorems with these explicit results.

\section{Critical points and the Gradient Flow}
Let us briefly recall the derivation of the critical points, see for example \cite{KohDiss} (Proposition 2.4).
\begin{Lem}[Critical points and Second Variation]
\begin{enumerate}
 \item The critical points of the energy functional (\ref{energy-functional}) satisfy
\begin{equation}
\label{first-variation}
\tau(\gamma)(s)=Z(\gamma')(s),
\end{equation}
where \(\tau\) is the tension field of the curve \(\gamma\) and \(Z\in\Gamma(\Hom(TN,TN))\).
\item The second variation of the energy functional (\ref{energy-functional}) yields
\begin{equation}
\label{second-variation}
\frac{\delta^2}{\delta\gamma^2}E(\gamma)=\int_{S^1}(|\nabla\eta|^2-\langle R^N(\eta,\gamma')\eta,\gamma'\rangle+\langle(\nabla_\eta Z)(\gamma'),\eta\rangle
+\langle Z(\nabla\eta),\eta\rangle)ds.
\end{equation}
\end{enumerate}
\end{Lem}
Solutions of (\(\ref{first-variation}\)) are called \emph{magnetic geodesics}.
\begin{proof}
Consider a family of smooth variations of $\gamma$ satisfying $\frac{\partial\gamma_t}{\partial t}\big|_{t=0}=\eta$.
The first variation of the energy of a curve is given by
\[
\frac{d}{dt}\bigg|_{t=0}\frac{1}{2}\int_{S^1}|\gamma_t'|^2ds=-\int_{S^1}\langle\tau(\gamma),\eta\rangle ds
\]
with variational vector field \(\eta\), see for example \cite{MR2431658}, p.2.

The first variation of the holonomy functional on the other hand gives 
\[
\frac{d}{dt}\bigg|_{t=0}\int_{S^1}\gamma_t^\ast A=-\int_{S^1}\Omega(\eta,\gamma')ds
\]
with variational vector field \(\eta\) and $\Omega = dA$, see for example \cite{MR2362847}, p.234.
To the 2-form \(\Omega\in\Gamma(\Lambda^2T^*N)\), we associate the smooth section $Z$ of \(\Hom(TN,TN)\) defined by the equation
\begin{equation} \label{EqZOmega}
\langle\eta,Z(\xi)\rangle=\Omega(\eta,\xi)
\end{equation}
for all \(\eta,\xi\in TN\) and thus the formula for the first variation follows.
\par\medskip
Concerning the second variation, remember the second variation of the energy of a curve, (see for example \cite{MR2431658}, p.8)
\[
\frac{d}{dt}\bigg|_{t=0}\frac{1}{2}\int_{S^1}|\gamma_t'|^2ds=\int_{S^1}(|\nabla\eta|^2-\langle R^N(\eta,\gamma')\eta,\gamma'\rangle)ds.
\]
For the second variation of the magnetic term, consider again a variation of \(\gamma\) satisfying 
$\frac{\partial\gamma_t}{\partial t}\big|_{t=0}=\eta$ and calculate
\[
\frac{\nabla}{\partial t}Z(\gamma'_t)=(\nabla_{\dot{\gamma_t}}Z)(\gamma_t')+Z(\nabla\dot{\gamma_t}).
\]
Evaluating at \(t=0\) yields the result.
\end{proof}
\begin{Bem}
If the magnetic field is not exact, then the first variation of \eqref{energy-functional-holonomy}
also gives \eqref{first-variation}.
\end{Bem}

Now we recall two standard Bochner-formulas, which were already proven in \cite{MR2551140}, Proposition 10.
\begin{Lem}[Bochner Formulas]
Let \(\gamma_t\colon S^1\times [0,T)\to N\) be a solution of (\ref{gradient-flow}) and let \(Z\in\Gamma(\Hom(TN,TN))\). 
Then the following Bochner formulas hold:
\begin{align}
\label{bochner1}
\frac{\partial}{\partial t}\frac{1}{2}|\gamma'_t|^2=&\Delta\frac{1}{2}|\gamma'_t|^2-|\tau(\gamma_t)|^2+\langle Z(\gamma'_t),\tau(\gamma_t)\rangle, \\
\label{bochner2}
\frac{\partial}{\partial t}\frac{1}{2}|\dot{\gamma}_t|^2=&\Delta\frac{1}{2}|\dot{\gamma}_t|^2-|\nabla\dot{\gamma}_t|^2+\langle R^N(\dot{\gamma}_t,\gamma'_t)\dot{\gamma}_t,\gamma'_t\rangle
-\langle(\nabla_{\dot{\gamma_t}}Z)(\gamma'_t),\dot{\gamma_t}\rangle-\langle Z(\frac{\nabla}{\partial t}\gamma'_t),\dot{\gamma}_t\rangle.
\end{align}
\end{Lem}

\begin{proof}
Regarding the first equation, a direct computation yields
\[
\frac{\partial}{\partial t}\frac{1}{2}|\gamma'_t|^2=\Delta\frac{1}{2}|\gamma'_t|^2-|\tau(\gamma_t)|^2-\langle \nabla Z(\gamma'_t),\gamma'_t\rangle.
\]
Using the identity
\[
\langle \nabla Z(\gamma'_t),\gamma'_t\rangle=\frac{\partial}{\partial s}\underbrace{\langle Z(\gamma'_t),\gamma'_t\rangle}_{=0}-\langle Z(\gamma'_t),\tau(\gamma_t)\rangle,
\]
the first assertion follows.
For the second statement, again by a direct computation we get
\[
\frac{\partial}{\partial t}\frac{1}{2}|\dot{\gamma}_t|^2=\Delta\frac{1}{2}|\dot{\gamma}_t|^2-|\nabla\dot{\gamma}_t|^2
+\langle R^N(\dot{\gamma}_t,\gamma'_t)\dot{\gamma}_t,\gamma'_t\rangle
-\langle\frac{\nabla}{\partial t} Z(\gamma'_t),\dot{\gamma}_t\rangle.
\]
Differentiating with respect to \(t\) we find
\[
\frac{\nabla}{\partial t}Z(\gamma'_t)=(\nabla_{\dot{\gamma_t}}Z)(\gamma_t')+Z(\frac{\nabla}{\partial t}\gamma_t')
\]
and thus the statement follows.
\end{proof}

With the help of the maximum principle we are now able to derive estimates,
which are obtained similar to Corollary 13 in \cite{MR2551140}:
\begin{Lem}
Let \(\gamma_t\colon S^1\times [0,T)\to N\) be a solution of (\ref{gradient-flow}) and 
let \(Z\in\Gamma (\Hom(TN,TN))\). Then the following estimates hold
\begin{align} 
\label{estimate-gamma'}
|\gamma'_t|^2\leq& |\gamma'_0|^2e^{C_1t},\\
\label{estimate-gamma-dot}
|\dot{\gamma}_t|^2\leq& |\dot{\gamma}_0|^2e^{C_2e^{t}+C_3t}.
\end{align}
The constant \(C_1\) depends on \(|Z|_{L^\infty}\), the constant \(C_2\) depends on \(N,|Z|_{L^\infty},|\nabla Z|_{L^\infty},|\gamma_0'|\)
and the constant \(C_3\) depends on \(|Z|_{L^\infty},|\gamma_0'|\).
\end{Lem}
\begin{proof}
To derive the first inequality we estimate the Bochner formula (\ref{bochner1}) and find
\[
\frac{\partial}{\partial t}\frac{1}{2}|\gamma'_t|^2\leq\Delta\frac{1}{2}|\gamma'_t|^2+\frac{1}{4}|Z|^2_{L^\infty}|\gamma'_t|^2.
\]
Now, apply the maximum principle and set \(c_1=\frac{1}{2}|Z|^2_{L^\infty}\).
Regarding the second inequality, we estimate the Bochner formula (\ref{bochner2}) and use the estimate on \(|\gamma_t'|^2\) to obtain
\begin{align}
\frac{\partial}{\partial t}\frac{1}{2}|\dot{\gamma}_t|^2 \leq& \Delta\frac{1}{2}|\dot{\gamma}_t|^2-|\nabla\dot{\gamma}_t|^2+c_2|\gamma'_t|^2|\dot{\gamma}_t|^2+c_3|\dot{\gamma}_t|^2|\gamma'_t|
+\sqrt{2c_1}|\dot{\gamma}_t||\frac{\nabla}{\partial t}\gamma'| \label{EqBochnerEstimate} \\
\leq&\Delta\frac{1}{2}|\dot{\gamma}_t|^2
+c_2|\dot{\gamma}_t|^2|\gamma'_0|^2e^{c_1t}+c_3|\dot{\gamma}_t|^2|\gamma'_0|e^{\frac{c_1}{2}t}+\frac{c_1}{2}|\dot{\gamma}_t|^2 \notag \\
\leq&\Delta\frac{1}{2}|\dot{\gamma}_t|^2
+|\dot{\gamma}_t|^2\left(c_2|\gamma'_0|^2e^{c_1t}+c_3|\gamma'_0|e^{\frac{c_1}{2}t}+\frac{c_1}{2}\right) \notag
\end{align}
with the constants \(c_2=|R^N|_{L^\infty}\) and \(c_3=|\nabla Z|_{L^\infty}\).
Again, by application of the maximum principle, we have to estimate the solution of
\[
\frac{\partial}{\partial t}|\dot{\gamma}_t|^2\leq c_4|\dot{\gamma}_t|^2(e^{c_1t}+c_5)
\]
with some positive constants \(c_4\) and \(c_5\). Integrating the ODE and rearranging the constants
completes the proof.
\end{proof}

The estimates from the last Lemma are sufficient to establish the following statement, which is proven in \cite{MR2551140}{, Theorem 1}:

\begin{Satz}\label{Dennislangzeit}
For any $\gamma_0\in C^{2+\alpha}(S^1,N)$, there exists a unique $\gamma \in C^{\infty}(S^1 \times [0,\infty),N)$ 
satisfying \eqref{gradient-flow} with initial condition $\gamma(\cdot,0)=\gamma_0$. 
\end{Satz}

\section{Convergence of the Gradient Flow} \label{SecConvergence}

If we want to achieve the convergence of the gradient flow, we have to establish energy estimates that are independent of the deformation
parameter \(t\). Thus, we have to improve the estimates obtained by the maximum principle \eqref{estimate-gamma'} and \eqref{estimate-gamma-dot}. 
We will make use of the following Lemma, which combines the pointwise maximum principle with an integral norm.

\begin{Lem}
\label{maximum-principle-l2}
Assume that \((M,h)\) is a compact Riemannian manifold. If a function \(u(x,t)\geq 0\) satisfies
\[
\frac{\partial u}{\partial t}\leq \Delta_h u+Cu,
\]
and if in addition we have the bound
\[
U(t)=\int_Mu(x,t)dM\leq U_0,
\]
then there exists a uniform bound on 
\[
u(x,t)\leq e^CKU_0 
\]
with the constant \(K\) depending only on the geometry of \(M\).
\end{Lem}
\begin{proof}
A proof can for example be found in \cite{MR2744149}, p.\ 284.
\end{proof}

Based on this result, we obtain the following general result concerning convergence of the heat flow:

\begin{Thm}[Convergence]
\label{Convergence}
Let \(\gamma:S^1\times [0,\infty)\to N\) be a solution of \eqref{gradient-flow} and \((N,g)\) be a compact Riemannian manifold.
If we have uniform bounds 
\begin{equation}
\int_{S^1}|\gamma'_t|^2ds\leq C,\qquad  \int_0^\infty\int_{S^1}|\dot{\gamma}_t|^2dsdt\leq C
\end{equation}
for all \(t\in [0,\infty)\), then the evolution equation \eqref{gradient-flow}
subconverges in \(C^2(S^1,N)\) to a magnetic geodesic \(\gamma_\infty\).
\end{Thm}

\begin{proof}
By assumption, we have a uniform bound on the \(L^2\) norm of \(\gamma'_t\).
Together with the pointwise equation (\ref{bochner1}) and Lemma \ref{maximum-principle-l2} we get the pointwise bound \(|\gamma'_t|^2\leq C\).
Using this bound on \(|\gamma'_t|^2\) together with an estimate analogous to \eqref{EqBochnerEstimate}, we note that \(|\dot{\gamma}_t|^2\) now satisfies 
\[
\frac{\partial}{\partial t}|\dot{\gamma}_t|^2\leq\Delta|\dot{\gamma}_t|^2+C|\dot{\gamma}_t|^2
\]
for some constant \(C\). Integrating over \(S^1\) and \(t\) from \(0\) to \(T\), we find the following bound
\begin{equation} \label{EqEstimate2}
\int_{S^1}|\dot{\gamma}_t|^2ds\leq \int_{S^1}|\dot{\gamma}_0|^2ds + \int_0^T\int_{S^1}|\dot{\gamma}_t|^2dsdt\leq C.
\end{equation}
By Lemma \ref{maximum-principle-l2} and the assumption, we now also obtain a pointwise bound on \(|\dot\gamma_t|^2\leq C\), uniform in $t$.\\
As a next step, we apply Nash's embedding theorem, that is we assume that $N$ is realized as a Riemannian submanifold of some $\R^q$. 
Thus, we may work with vector valued maps $S^1 \times [0,\infty) \rightarrow \R^q$ and $S^1 \rightarrow \R^q$ whose image is contained in $N \subset \R^q$ 
and use the associated H\"older spaces $C^{k,\alpha}(S^1,\R^q)$. Equation \eqref{gradient-flow} now takes the form 
\begin{align} \label{EqGradientFlowRn}
\dot{\gamma} &= -\Delta_{S^1,\R^n}(\gamma) + \mathbb{I}^N(\gamma',\gamma') + Z(\gamma'), 
\end{align}
where $\Delta_{S^1,\R^n}$ denotes the linear Hodge-Laplacian acting on $C^\infty(S^1,\R^q)$ and $\mathbb{I}^N$ the second fundamental form of the embedding $N \subset \R$. 
We improve our estimates with the help of Schauder theory. Following the presentation in \cite{MR2551140}{(proof of Theorem 22, in particular (31))} and viewing \eqref{EqGradientFlowRn} as a linear elliptic equation for $\gamma$, 
the elliptic Schauder estimates (\cite{MR2371700}, p.463, Thm.27) in $\R^q$ yield 
\begin{align} \label{Eq_1_alpha_Reg}
 |\gamma_t|_{C^{1,\alpha}(S^1,\R^q)} &\leq C \bigl( |\gamma_t'|^2_{L^\infty(S^1,\R^q)} + |\gamma_t'|_{L^\infty(S^1,\R^q)} + |\gamma_t|_{L^\infty(S^1,\R^q)} + |\dot{\gamma}_t|_{L^\infty(S^1,\R^q)} \bigr),
\end{align}
where the constant $C$ may depend on $N,Z, \gamma_0, \alpha$ and the embedding $N \hookrightarrow \R^q$ but not on $t \in [0,\infty)$. By the first part of the proof and the compactness of $N$, we thus get the bound
\begin{align*}
 |\gamma_t|_{C^{1,\alpha}(S^1,\R^q)} &\leq C  
\end{align*}
uniform in $t$. Viewing \eqref{EqGradientFlowRn} as a linear parabolic equation for $\gamma$ and using the corresponding Schauder estimates 
(see again \cite{MR2551140}, proof of Theorem 22 and also \cite{MR1465184}, Theorem 4.9 for the local version of the estimate), we obtain
\begin{align*}
 |\gamma_t|_{C^{2,\alpha}(S^1,\R^q)} + |\dot{\gamma}_t|_{C^\alpha(S^1,\R^q)} 
&\leq C \bigl( |\mathbb{I}^N(\gamma_t',\gamma_t') + Z(\gamma_t')|_{C^\alpha(S^1,\R^q)} + |\gamma_t|_{L^\infty(S^1,\R^q)} \bigr) \\
&\leq C \bigl( |\gamma_t'|^2_{C^{\alpha}(S^1,\R^q)} + |\gamma_t'|_{C^\alpha(S^1,\R^q)} + |\gamma_t|_{L^\infty(S^1,\R^q)} \bigr).
\end{align*}
Using the compactness of $N$ and \eqref{Eq_1_alpha_Reg}, we get bounds
\begin{align} \label{EqUniformBound1}
|\gamma_t|_{C^{2,\alpha}(S^1,\R^q)}, |\dot{\gamma}_t|_{C^\alpha(S^1,\R^q)} &\leq C ,
\end{align}
which are again uniform in $t \in [0,\infty)$. Now, by assumption, we have the estimate 
\begin{align*} 
\int_0^{\infty}\int_{S^1}|\dot{\gamma_t}|^2dsdt\leq C.
\end{align*}
Hence, there exists a sequence $t_k \rightarrow \infty$ such that 
\begin{align} \label{EqL2Convergence} 
|\dot{\gamma}_{t_k}|^2_{L^2(S^1)}\to 0 
\end{align}
as $k \rightarrow \infty$. Using the bounds from \eqref{EqUniformBound1} and the Theorem of Arzela and Ascoli, there exists a subsequence (again denoted by $\gamma_{t_k}$), 
which converges in $C^2$ to a limiting map $\gamma_\infty \in C^{2}(S^1,\R^q)$ such that $\dot{\gamma}_{t_k}$ converges in $C^0(S^1,\R^q)$. 
In fact, $\gamma_\infty$ defines an element of $C^2(S^1,N)$ because we have $\gamma_t(S^1) \subset N \subset \R^q$ for all $t$. 
Finally, \eqref{EqL2Convergence} implies that $\lim_{k \to \infty}\dot{\gamma}_{t_k} = 0$ and we conclude that $\gamma_\infty$ is a $C^2$-solution of \eqref{first-variation}. 
Since $\gamma$ is smooth in $t$, it is moreover clear that $\gamma_\infty$ is homotopic to $\gamma_0$.
\end{proof}

Since the limit in the previous proof is a $C^2$-solution of \eqref{first-variation} by construction
and recalling that a magnetic geodesic has constant energy, which follows from
\[
\frac{\partial}{\partial s}\frac{1}{2}|\gamma'_\infty|^2=\langle\tau(\gamma_\infty),\gamma'_\infty\rangle=\langle Z(\gamma'_\infty),\gamma'_\infty\rangle=0,
\]
the standard bootstrap argument for elliptic equations now implies the following regularity result:

\begin{Cor}
The limit $\gamma_\infty$ from Theorem \ref{Convergence} is smooth.
\end{Cor}

\begin{Bem}
In the case of the heat flow for \emph{geodesics} and under the assumption that the target manifold \(N\) has non-positive curvature,
a theorem due to Hartmann \cite{MR0214004} states that the limiting map \(\gamma_\infty\) is independent of
the chosen subsequence. This theorem uses the fact, that the second 
variation of the harmonic energy is positive. Due to the extra term in the energy
functional \eqref{energy-functional}, there is a corresponding term in the second variation \eqref{second-variation} and we do not get such a statement here.
\end{Bem}

\begin{Bem}
So far, we have established the existence of a convergent subsequence.
This does not necessarily imply that the gradient flow converges itself.
This phenomena also occurs in the heat flow for closed geodesics,
see \cite{choi-parker}.
\end{Bem}

\begin{Bem} \label{RemGeneralConvergence}
Let \(\gamma_t\) be a solution of \eqref{gradient-flow}, then for $T > 0$, we have 
\begin{equation}
\label{energy-equality}
\frac{1}{2}\int_{S^1}|\gamma_T'|^2ds+\int_0^T\int_{S^1}|\dot{\gamma}_t|^2dsdt+\int_0^T\int_{S^1}\Omega(\dot{\gamma}_t,\gamma_t')dsdt=\frac{1}{2}\int_{S^1}|\gamma'_0|^2ds.
\end{equation}
The only term in \eqref{energy-equality} not having a definite sign is the 
last term on the left hand side. By Theorem \ref{Convergence}, we may in particular expect that if we have a uniform bound 
\begin{equation*}
\int_0^T\int_{S^1}\Omega(\dot{\gamma}_t,\gamma_t')dsdt\geq C > -\infty,
\end{equation*}
then the evolution equation \eqref{gradient-flow} subconverges to a magnetic geodesic.
\end{Bem}

In the following we will describe two situations in which we can ensure
the necessary estimates such that we can apply Theorem \ref{Convergence}.

\subsection{The case of null-homotopic curves}
Throughout this section we assume that the curves \(\gamma_t\) are null-homotopic. If \(U\) is a subset of \(N\), we define
\begin{align} \label{EqDefKappa_r}
\kappa_U &:=\sup\{\langle R^N(X,Y)Y,X\rangle\mid X,Y\in T_pM,|X|=|Y|=1,p\in U\}
\end{align}

To give a sufficient condition for the application of Theorem \ref{Convergence}, we want to apply the following Poincaré type inequality (for a derivation see \cite{MR834094}, p.58):
\begin{Satz}
\label{ThmOttarson}
Let \(\gamma\colon S^1\to N\) be a non-constant closed \(C^2\) curve and assume
w.l.o.g. that the maximum of the energy density occurs at \(s=0\). Suppose that the image
of \(\gamma\) is contained in a ball \(B_r(\gamma(0))\) such that the exponential map 
\[
\exp_{\gamma(0)}\colon T_{\gamma(0)}N\to N, 
\]
restricted to \(B_r(0)\) in \(T_{\gamma(0)}N\), is a diffeomorphism onto \(B_r(\gamma(0))\). In case \(\kappa_B>0\), assume in addition that \(r<(2\sqrt{\kappa_B})^{-1}\).
Then the following inequality holds:
\begin{equation}
\label{poincare-tension}
\frac{1}{4\pi^2}\int_{S^1}|\gamma'|^2ds\leq\int_{S^1}|\tau(\gamma)|^2ds
\end{equation}
\end{Satz}

\begin{Bem} \label{RemOttarson}
\begin{enumerate}
 \item The estimate in \eqref{poincare-tension} can be improved for specific geometries. As an example, consider $\R^n$ with the flat Euclidean metric. 
 A simple calculation using Fourier series (see also \cite{MR834094} p.57) then shows that the inequality still holds with $1/4\pi^2$ replaced by $1$. 
 For a general Riemmanian manifold $(N,g)$, let us denote the optimal constant for \eqref{poincare-tension} by $O(N,g)$. By Theorem \ref{ThmOttarson}, we thus have $O(N,g) \in [1/(4\pi^2), \infty)$. 
 Moreover, computing the	 integrals for $\gamma$ a circle with radius 1, it is straightforward to see that actually $O(\R^n,g_{Eucl}) = O(T^n,g_{flat}) = 1$.
 \item The following examples show that the assumptions in Theorem \ref{ThmOttarson} are necessary: Consider the two-dimensional flat torus and a nontrivial closed geodesic on it. 
 Since \eqref{poincare-tension} is clearly violated, we see that the condition concerning the exponential map cannot be omitted. 
 Similarly, the statement is not valid without the bound on the radius $r$ in case of $\kappa_{B_r} > 0$: Consider a circle $c_h$ of latitude $h > 0$ on $S^2$. 
 Clearly, for each point $p$ on $c_h$, the curve is completely contained in the ball of radius $\pi$ (= the injectivity radius of $S^2$) around $p$.
 On the other hand, for $h \searrow 0$, this circle converges to the equator and hence $\int_{S^1} |\tau(c_h)|^2ds \rightarrow 0$, whereas $\int_{S^1} |c_h'|^2ds$ converges to some positive number. 
 This shows that we have to restrict to balls of radius smaller than the injectivity radius in order to ensure \eqref{poincare-tension} without assuming $K \leq 0$. 
\end{enumerate}
\end{Bem}

\begin{Lem} \label{LemZAssump}
Let \(\gamma:S^1\times [0,\infty)\to N\) be a solution of (\ref{gradient-flow}) and \((N,g)\) be a compact Riemannian manifold.
Assume that
\begin{equation}
|Z|^2_{L^\infty}\leq\frac{1}{4\pi^2}
\end{equation}
and that Theorem \ref{ThmOttarson} can be applied to the curves $\gamma_t$ for all $t \in [0,\infty)$.
Then we obtain the uniform bounds  
\begin{equation}
\int_{S^1}|\gamma'_t|^2ds\leq C,\qquad  \int_0^\infty\int_{S^1}|\dot{\gamma}_t|^2dsdt\leq C.
\end{equation}
\end{Lem}
\begin{proof}
We rewrite (\ref{energy-equality}) in the following way
\begin{align*}
\frac{1}{2}\int_{S^1}|\gamma'_T|^2ds+\frac{1}{2}\int_0^T\int_{S^1}|\dot{\gamma}_t|^2dsdt=
\frac{1}{2}\int_{S^1}|\gamma'_0|^2ds-\frac{1}{2}\int_0^T\int_{S^1}|\dot{\gamma}_t|^2dsdt-\int_0^T\int_{S^1}\Omega(\dot{\gamma}_t,\gamma_t')dsdt.
\end{align*}
By a direct calculation using the evolution equation (\ref{gradient-flow}) we obtain
\[
-\frac{1}{2}|\dot{\gamma}_t|^2-\Omega(\dot{\gamma}_t,\gamma'_t)=-\frac{1}{2}|\tau(\gamma_t)|^2+\frac{1}{2}|Z(\gamma_t')|^2.
\]
Combining both equations, we find for \(T>0\)
\[
\int_{S^1}|\gamma'_T|^2ds+\int_0^T\int_{S^1}|\dot{\gamma}_t|^2dsdt=\int_{S^1}|\gamma'_0|^2ds + \int_0^T\int_{S^1}(-|\tau(\gamma_t)|^2+|Z(\gamma'_t)|^2)dsdt.
\]
Using the assumption on \(|Z|^2_{L^\infty}\) and applying Ottarsson's inequality \eqref{poincare-tension}, we obtain
\begin{equation}
\int_{S^1}|Z(\gamma'_t)|^2ds-\int_{S^1}|\tau(\gamma_t)|^2ds\leq(4\pi^2|Z|^2_{L^\infty}-1)\int_{S^1}|\tau(\gamma_t)|^2ds\leq 0,
\end{equation}
which gives the result.
\end{proof}

\begin{Bem}
If we think of magnetic geodesics as curves of prescribed geodesic curvature, then the assumption of Lemma \eqref{LemZAssump}
means that we have to restrict to small curvature.
\end{Bem}

We next discuss the question whether the conditions from Theorem \ref{ThmOttarson} are preserved under the flow, i.e. whether we may use estimate \eqref{poincare-tension} provided it is satisfied for $t =0$. 
The following Lemma shows, that this is indeed possible for curves of sufficiently low kinetic energy. To give a precise statement of this condition, let us define 
\begin{align}
r(N) &:= \begin{cases}
          \mathrm{injrad}(N) &\text{ if } \kappa_N \leq 0, \\
          \mathrm{min}(\mathrm{injrad}(N),(2\sqrt{\kappa_N })^{-1} ) &\text{ if } \kappa_N > 0,  
         \end{cases}   \label{EqDefRadius}
\end{align}
where $\kappa_N$ was defined in \eqref{EqDefKappa_r}.

\begin{Lem} \label{LemApplyOttarson}
If the assumptions from Lemma \ref{LemZAssump} on $Z$ are satisfied, then for any initial curve of kinetic energy smaller than $r(N)^2 / (16\pi)$, 
Ottarson's estimate \eqref{poincare-tension} is valid for all curves $\gamma_t, \ t \in [0,\infty)$ of the associated flow.
\end{Lem}

\begin{proof}
We first observe that  
\begin{align}
\frac{\partial}{\partial t}\frac{1}{2}\int_{S^1}|\gamma'_t|^2ds=&-\int_{S^1}|\tau(\gamma_t)|^2ds+\int_{S^1}\langle\tau(\gamma_t),Z(\gamma'_t)\rangle ds \notag \\
\leq&\frac{1}{2}\big(-\int_{S^1}|\tau(\gamma_t)|^2ds+|Z|^2_{L^\infty}\int_{S^1}|\gamma'_t|^2ds\big) \notag \\
\leq&\frac{1}{2}\big(|Z|_{L^\infty}^2-\frac{1}{4\pi^2}\big)\int_{S^1}|\gamma'_t|^2ds, \label{EqEstimate1}
\end{align}
where we applied (\ref{poincare-tension}) in the last step. By the definition of $r(N)$ from \eqref{EqDefRadius} and the compactness of $N$, we have $r(N) > 0$. 
Let $\gamma_0$ be an initial curve satisfying $E_{kin}(\gamma_0) < (r(N))^2 / (16\pi)$. 
By the Cauchy-Schwarz inequality, the length of $\gamma_0$ can be estimated by $L(\gamma_0) < r(N)/2$. 
Thus, for any point $p \in \gamma_0(S^1)$, we have $\gamma_0(S^1) \subset B_{r(N)}(p)$ and by the choice of $r(N)$ and continuity, 
we may apply Theorem \ref{ThmOttarson} on some time interval $[0,\epsilon)$. 
Now assume that 
\begin{align*}
J &:= \{t \in [0,\infty) \mid \eqref{poincare-tension} \text{ is not valid } \}  \subset [0,\infty)
\end{align*}
is \emph{not empty}. Then, $T := \inf(J) \geq \epsilon > 0$ and \eqref{poincare-tension} holds on $[0,T)$. 
By \eqref{EqEstimate1} and the assumption on $Z$, we have $\partial E_{kin}(\gamma_t) / \partial t \leq 0$ on $[0,T)$ which implies $E(\gamma_T) \leq E(\gamma_0)$. 
Thus, we may argue as before to show that \eqref{poincare-tension} is in fact valid on $[0, T+\epsilon')$. 
This contradiction proves that $J$ must be empty, i.e. \eqref{poincare-tension} holds for all $t > 0$.
\end{proof}

Combining the results of Theorem \ref{Convergence} and Lemmas \ref{LemZAssump}, \ref{LemApplyOttarson}, we have shown

\begin{Satz} \label{ThmConvergenceOttarson}
Let $Z$ satisfy $|Z|^2_{L^\infty} \leq 1/(4\pi^2)$ and assume that $E_{kin}(\gamma_0) < r(N)^2/(16\pi)$. 
Then the solution of \eqref{gradient-flow} subconverges to a magnetic geodesic. 
\end{Satz}
Since all curves satisfying the assumptions of Theorem \ref{ThmOttarson} are contractible, the question of nontriviality of the limit curve $\gamma_\infty$ arises naturally. 
To discuss this question, we first observe that integrating the estimate \eqref{EqEstimate1} with respect to $t$ directly yields the following estimate on the kinetic energy under the flow: 

\begin{Lem} \label{LemEkinEstimate}
Let \(\gamma:S^1\times [0,\infty)\to N\) be a solution of (\ref{gradient-flow}) and \((N,g)\) be a compact Riemannian manifold.
If Theorem \ref{ThmOttarson} can be applied to the curves $\gamma_t$ for all $t \in [0,\infty)$, then we have
\begin{equation}
\int_{S^1}|\gamma'_t|^2ds\leq e^{\big(|Z|_{L^\infty}^2-\frac{1}{4\pi^2}\big)t}\int_{S^1}|\gamma'_0|^2ds.
\end{equation}
\end{Lem}

This allows us to conclude that the flow in fact converges to a point in many cases:
\begin{Cor} \label{trivial-limit}
If the magnetic field $Z$ satisfies
\begin{equation}
|Z|^2_{L^\infty} < \frac{1}{4\pi^2},
\end{equation}
then the flow converges to the trivial magnetic geodesic \(\gamma_\infty\). 
In fact, under this assumption on $Z$, any magnetic geodesic that satisfies \eqref{poincare-tension} is trivial.
\end{Cor}

\begin{proof}
Lemma \ref{LemEkinEstimate} directly implies, that $E_{kin}(\gamma_\infty)=0$, i.e. $\gamma_\infty$ is constant.
An argument analogous to the one used in \cite{choi-parker} (Chapter 9) shows that under this condition, the flow does not only subconverge but actually converges. 
Moreover, if a curve $\gamma$ satisfies \eqref{first-variation} together with \eqref{poincare-tension}, we have $\int_{S^1} |\tau(\gamma)|^2ds \leq \frac{1}{4\pi^2}|Z|^2_{L^\infty} \int_{S^1} |\tau(\gamma)|^2ds$.
Hence, $\int|\tau(\gamma)|^2 ds= 0$ and \eqref{poincare-tension} implies that $\gamma$ is constant.
\end{proof}

\begin{Bem}
The situation described in Corollary \ref{trivial-limit} is similar to the one for the heat flow for geodesics, see \cite{MR834094}, Theorem 4A and also \cite{choi-parker}, Section 9.
\end{Bem}

\begin{Bem}\label{RemZ_Eins_Optimal}
As noted in Remark \ref{RemOttarson} (1), \eqref{poincare-tension} may be improved using the optimal constant $O(N,g)$. 
From the argument given in Corollary \ref{trivial-limit}, we see that we obtain convergence to a trivial limit in case $|Z|^2_{L^\infty} < O(N,g)$. 
A non-trivial limit $\gamma_\infty$ can only be obtained, if $|Z|^2_{L^\infty} = O(N,g)$ (note that for larger values, Theorem \ref{ThmConvergenceOttarson} may no longer be applied!). 
In general, it is difficult to say what happens in this case of equality. However, we will later see from Example \ref{ExTorusPart1} that our result is in fact optimal in the following sense: 
In accordance with Corollary \ref{trivial-limit}, the flow converges to a point for $|Z|^2_{L^\infty} < 1 = O(T^2, g_{flat})$. 
For $|Z|^2_{L^\infty} = 1$, we do in fact obtain a nontrivial limit and for $|Z|_{L^\infty} > 1$, the flow does not converge in general. 
Finally, all periodic magnetic geodesics in this example are in fact contractible.
\end{Bem}

Theorem \ref{Dennislangzeit} and Theorem \ref{ThmConvergenceOttarson} in fact also allow us to conclude the existence of nontrivial periodic magnetic geodesics via the heat flow method without assuming that $Z$ is small.
Let $\Lambda > 0$ and $\gamma : S^1 \times [0,\infty) \rightarrow N$ a solution of \eqref{gradient-flow}, defined on the circle of length $2\pi$. 
Then, a straightforward calculation shows that the rescaled curve $\gamma_\Lambda : S^1_{\Lambda^{-1}} \times [0,\infty), \gamma_\Lambda(s,t) := \gamma(\Lambda s, \Lambda^2 t)$ satisfies the rescaled equation
\begin{align} \label{rescaledEOM}
 \dot{\gamma}_\Lambda &= \tau(\gamma_\Lambda) - \Lambda Z(\gamma_\Lambda).
\end{align}

Using this rescaling technique, we obtain 

\begin{Cor} \label{CorRescaling}
Given an arbitrary $Z \in \Gamma(\Hom(TN,TN))$, there exist solutions to \eqref{gradient-flow} which converge to (possibly trivial) closed magnetic geodesics with respect to $Z$, defined on a possibly rescaled loop.
\end{Cor}

\begin{proof}
Let $O(N,g)$ denote the constant from Remark \ref{RemOttarson} (1), choose $\Lambda \leq O(N,g)/|Z|^2_{L^\infty}$ and set $Z_\Lambda := \Lambda Z$ .
Then Theorem \ref{Dennislangzeit} guarantees a solution $\gamma_\Lambda : S^1 \times [0,\infty)$ to \eqref{rescaledEOM}. 
By Theorem \ref{ThmConvergenceOttarson} and the choice of $\Lambda$, this flow subconverges (as $t \rightarrow \infty$) to $\gamma_{\Lambda,\infty}$, 
which  satisfies $\tau(\gamma_{\Lambda,\infty}) = Z_{\Lambda}(\gamma'_{\Lambda,\infty})$.
Scaling by $\Lambda^{-1}$ as described above, the resulting map $\gamma : S^1_{\Lambda} \times [0,\infty), \gamma(s,t) := \gamma_\Lambda(\Lambda^{-1} s, \Lambda^{-2} t)$ solves the original equation \eqref{gradient-flow}. 
It is clear that $\gamma$ still subconverges to a solution of the equation for magnetic geodesics with the original $Z$ which, however, is defined on $S^1_{\Lambda}$.
\end{proof}

Of course, choosing $\Lambda < O(N,g)/|Z|^2_{L^\infty}$, Corollary \ref{trivial-limit} implies, that the limit will aways be trivial.

\begin{Bem}
It is well know that in our (approximately flat) physical space $\R^3$, periodic magnetic geodesics - in fact circles - can be observed for \emph{homogeneous} magnetic fields of (almost) arbitrary strength $B_0$. 
This is included in Corollary \ref{CorRescaling} but not in Theorem \ref{ThmConvergenceOttarson} if $B_0$ is large. 
In fact, it is easy to see that we have to perform a rescaling to include these solutions: 
Setting charge and mass to $1$, a straightforward calculation shows that equality of Lorentz- and centrifugal force translates into $2\pi/T = B_0$ where $T$ denotes the time of revolution. 
Since $T$ corresponds to the length of $S^1_\Lambda$, it is obvious that its length has to be ``adjusted''.  
\end{Bem}

\subsection{The case of an exact magnetic field}
In this section, we discuss the convergence of the gradient flow
in the case that the magnetic field is exact. Under this assumption we 
can exploit the fact that \eqref{gradient-flow} is derived from a the well defined energy \(E(\gamma)\).

\begin{Satz} \label{ThmConvExactCase}
If the two-form $\Omega$ is exact, then 
the solution $\gamma$ of \eqref{gradient-flow} subconverges to a magnetic geodesic \(\gamma_\infty\).
\end{Satz}

\begin{proof}
Using the fact that our evolution equation \eqref{gradient-flow}
is the \(L^2\)-gradient flow of \(E(\gamma)\) we obtain
\begin{equation}
\label{energy-inequality-exact}
E(\gamma_T)+\int_0^T\int_{S^1}|\dot\gamma_t|^2dsdt=\frac{1}{2}\int_{S^1}|\gamma'_T|^2ds+\int_{S^1}\gamma_T^\ast Ads+\int_0^T\int_{S^1}|\dot\gamma_t|^2dsdt=E(\gamma_0).
\end{equation}
Thus, we may directly conclude that 
\begin{equation}
\label{exact-energy-inequality}
\frac{1}{4}\int_{S^1}|\gamma'_T|^2ds+\int_0^T\int_{S^1}|\dot\gamma_t|^2dsdt\leq E(\gamma_0)+\vol(S^1)|A|^2_{L^\infty}.
\end{equation}
Using \eqref{bochner1} and \eqref{bochner2} together with \eqref{maximum-principle-l2}, we obtain a uniform bound on \(|\gamma_t'|^2\) and \(|\dot\gamma_t|^2\) from \eqref{exact-energy-inequality}.
Thus, the assumptions of Theorem \ref{Convergence} are satisfied and the evolution equation \eqref{gradient-flow} subconverges in \(C^2(S^1,N)\) to a magnetic geodesic \(\gamma_\infty\). 
Since $\gamma$ is smooth, it is clearly homotopic to $\gamma_0$. 
\end{proof}

Note that we did not use Ottarson's Theorem \ref{ThmOttarson} and henceforth do not have to restrict to contractible curves.
Starting with initial conditions in a prescribed homotopy class, we thus obtain 

\begin{Cor}
If $\Omega$ is exact, then each homotopy class of curves in $N$ contains a magnetic geodesic.
\end{Cor}

Note that in contrast to the heat flow for geodesics on manifolds with negative curvature, the magnetic geodesic need not be unique in its homotopy class. 
The argument which is used for the heat flow of ordinary geodesics (see \cite{MR0214004}) is no longer valid because of the presence of the magnetic contribution to the energy.
We will see an explicit counterexample in Example \ref{ExHyperbolic}: $\mathbb{H}^2$ (or compact quotients thereof), equipped with a suitable magnetic field $Z$, allows for contractible nontrivial closed magnetic geodesics.
The following statement shows that under additional assumptions on the initial conditions, the limit $\gamma_\infty$ will be nontrivial, even in the contractible case.

\begin{Prop} \label{LemNonTriv1}
Let $\Omega= dA$ be exact and  $\gamma_\infty$ be the limit of the associated flow (which exists by Theorem \ref{ThmConvExactCase}) with initial condition $\gamma_0$. 
Moreover, assume $E(\gamma_0) \leq 0$ and in addition $\gamma_0 \neq \mathrm{constant \ curve}$ if $E(\gamma_0)=0$. 
Then, $\gamma_\infty$ is a nontrivial closed magnetic geodesic.  
\end{Prop}

\begin{proof}
Let $\gamma$ be the solution to \eqref{gradient-flow}. By \eqref{energy-inequality-exact} applied to the interval $t \in [t_1,t_2]$, we have $E(\gamma_{t_2})\leq E(\gamma_{t_1})$ for any $t_2 \geq t_1 $ 
and hence, $E(\gamma_\infty) \leq E(\gamma_t)$ for any $t \geq 0$. If $E(\gamma_t) < 0$ for some $t$ (in particular for $t=0$), this proves the claim since the energy of a constant curve clearly vanishes. 
If $E(\gamma_t) = 0$ for all $t$, then we conclude from \eqref{energy-inequality-exact} that $\dot{\gamma}=0$ and $\gamma_\infty = \gamma_0$, which is nontrivial by assumption in this case. 
\end{proof}

The existence of contractible initial conditions $\gamma_0$ satisfying $E(\gamma_0) \leq 0$ can in fact be ensured provided the field does not vanish:

\begin{Cor}
Assume that $\Omega = dA \neq 0$. Then there exists $\gamma_0 : S^1_\Lambda \rightarrow N$ for some $\Lambda > 0$ such that the resulting flow from Theorem \ref{Dennislangzeit} subconverges to a nontrivial magnetic geodesic.  
\end{Cor}

\begin{proof}
By assumption, $A \in \Omega^1(N)$ for $\Omega$ cannot be closed. Thus, we can always find a small ball $B \subset N$ and a loop $c : S^1 \rightarrow B$ s.t. $\int_{S^1} c^\ast A \neq 0$. 
Reversing the orientation of $c$ if necessary, we may assume $\int_{S^1} c^\ast A < 0$. Rescaling $c$ as described above \eqref{rescaledEOM},
we see that $E_{kin}(c_\Lambda)$ becomes arbitrary small if $\Lambda$ tends to zero, whereas the term $\int_{S^1_\Lambda} (c_\Lambda)^\ast A$ is independent of $\Lambda$. 
Putting $\gamma_0 := c_\Lambda$ for an appropriate choice of $\Lambda$ yields $E(\gamma_0) \leq 0$ and we can apply Proposition \ref{LemNonTriv1}.
\end{proof}

\begin{Bem}
\begin{enumerate}
 \item The previous result is similar to the one stated in Theorem 1 of \cite{MR2679767}, which however uses a slightly different functional.
 \item In case that $\Omega$ is not exact, the argument used in the proof of Proposition \ref{LemNonTriv1} breaks down since there is no longer a well defined energy. 
 Of course, $\Omega$ always admits local potentials $A_{loc}$ and the argument clearly remains valid provided one can show that the curves $\gamma_t$ stay inside the domain of $A_{loc}$. 
 The examples in Section \ref{SecExamples} show that this may or may not be the case, see Remark \ref{RemTorusEx} \eqref{RemTorusExEnergy} and Remark \ref{RemSphere}.
 We are currently not aware of a condition ensuring the applicability of Proposition \ref{LemNonTriv1} in the non-exact case.
\end{enumerate}
\end{Bem}

\section{Examples} \label{SecExamples}
In this section we discuss some explicit examples on the two-dimensional torus $T^2$, on the two-dimensional sphere \(S^2\), and in the hyperbolic plane \(\mathbb{H}^2\). 
They illustrate the influence of different background geometries and initial conditions on the heat flow. 
Moreover, we believe it is helpful to discuss the abstract results from Section \ref{SecConvergence} in view of these examples.

\subsection{Magnetic geodesics on the torus} \label{ExTorus}
A general existence result for magnetic geo\-de\-sics on the two-dimensional flat torus has been obtained using methods from symplectic geometry in \cite{MR2036336} (Section 5). 
The heat flow in this specific case was studied in \cite{MR2551140}, p.457ff. Here, we want to make some further remarks concerning this example.\\

Let $C := S^1 \times \R \subset \R^3$ denote the cylinder of radius $1$. Using cylindrical coordinates $(\varphi,z) \in  (-\pi,\pi) \times \R$, a curve on $C$ depending on some parameter $t$ may be represented as 
\begin{align*}
\gamma_t &= (\cos(\varphi(s,t)),\sin(\varphi(s,t)),z(s,t)).
\end{align*}
The magnetic field $Z$ at $x \in C$ is defined by taking the vector product (in $\R^3$) of $v \in T_xC$ and the vector $B(x)$:
\begin{align} \label{EqTorusZ}
Z_x(v) &:= v \times B(x), & B=B_0(\cos(\varphi(x),\sin(\varphi(x)),0).
\end{align}
Here, $B_0 \in \R$ is a constant describing the strength of the magnetic field. Since $B(x) \perp T_xC$, we clearly have $Z_x \in \mathrm{End}(T_xC)$. 
It follows directly from \eqref{EqTorusZ} that $\nabla^C Z = 0$ and thus, the associated two-form $\Omega$ (cf. \ref{EqZOmega}) is closed. 
Finally, taking the quotient by $\Z$ in z-direction, we obtain a field $Z \in \mathrm{End}(TT^2)$ enjoying the same properties. 
Note that the resulting induced field strength in $\Omega^2(T^2)$ is given by $\Omega_{T^2} = B_0 vol_{T^2}$. In particular, it is not exact.
On the opposite, it is easy to check that on $C$, we have $\Omega_{C} = -d(z d\varphi)$, i.e. there exists a global potential $A = -z d\varphi \in \Omega^1(C)$.\\

Before solving the flow equation, we briefly discuss the solution to the ODE \eqref{first-variation} for magnetic geodesics in case $B_0 \neq 0$. 
This is most easily done on the level of the universal coverings $\R^2 \rightarrow T^2$ and $\R^2 \rightarrow C$, respectively. 
In fact, viewing $\R^2$ as subspace of $\R^3$, we have to determine the trajectories subject to a magnetic field of constant strength $B_0$ which is perpendicular to $\R^2$. 
It is well known that all the trajectories are circles in $\R^2$ of radius $|B_0|$, so in particular they are closed and contractible. 
Projecting  to $T^2$ or $C$, we conclude that all magnetic geodesics on these spaces subject to the magnetic field from \eqref{EqTorusZ} are closed and contractible.\\ 

The heat flow equation \eqref{gradient-flow} is equivalent to the following system of partial differential equations
\[
\dot\varphi=\varphi''-B_0z',\qquad \dot z=z''+B_0\varphi',
\]
which can be rewritten in terms of the complex variable \(\xi=\varphi+iz\) as
\begin{equation} \label{evolution-xi}
\dot{\xi}=\xi''+iB_0 \xi'.
\end{equation}
The integrand of the magnetic term in the energy identity \eqref{energy-equality}, which determines the asymptotic behaviour of the system, can be computed explicitly:
\[
\Omega(\dot{\gamma},\gamma')=\langle \dot{\gamma},Z(\gamma')\rangle=B_0\langle \dot{\gamma},\gamma'\times B\rangle=B_0(z'\dot{\varphi}-\dot{z}\varphi').
\]
We can integrate \eqref{evolution-xi} directly and obtain a family of solutions 
\begin{equation} \label{Eq_k_Mode}
\xi_k=e^{iks}e^{-(B_0k+k^2)t},
\end{equation}
where $k\in\N$. The most general solution is obtained by summing over all \(k\). Similar to \cite{MR2551140} we now solve this system for different types of initial data.\\

\begin{enumerate}
 \item \label{ExTorusPart1} Let $k \in \Z\setminus\{0\}$. Prescribing initial conditions $(\varphi_k(s,0),z_k(s,0))=(a\cos (ks),b\sin (ks))$ for $a,b> 0$, we find the following solution:
\begin{align*}
\varphi_k(s,t)&=\frac{a+b}{2}\cos (ks) e^{-(kB_0+k^2)t}+\frac{a-b}{2}\cos (ks) e^{-(-kB_0+k^2)t},\\
z_k(s,t)&=\frac{a+b}{2}\sin (ks)e^{-(kB_0+k^2)t}-\frac{a-b}{2}\sin (ks) e^{-(-kB_0+k^2)t}.
\end{align*}
Putting $B_0=k=1$, we reproduce the solution in \cite{MR2551140}, Example 8a) and it is easy to see that the flow converges as $t \to \infty$. 
More general, a brief inspection of the exponents shows that the flow (sub)-converges if and only if $|B_0| \leq |k|$. 
In case $|B_0| < |k|$, the flow in fact shrinks to a point and we only obtain a nontrivial limit for $B_0 = \pm k$. 
It is not hard to check that these limits are in fact magnetic geodesics. 
Explicit computation of the magnetic term in the energy identity \eqref{energy-equality} yields
\begin{align*}
\int_0^T\int_{S^1}\Omega(\dot{\gamma},\gamma')dsdt = -\pi abkB_0 
&- \tfrac{\pi}{4}(b-a)^2kB_0 e^{-2T(k^2-kB_0)} \\
&+ \tfrac{\pi}{4}(b+a)^2kB_0 e^{-2T(k^2+kB_0)}.
\end{align*}
We see that in accordance with Remark \ref{RemGeneralConvergence}, this term is bounded below with respect to $T$ if and only if $|B_0|\leq k$, i.e. if and only if the flow (sub)-converges. 

\item \label{ExTorusPart2} Prescribing the initial conditions $(\varphi(s,0),z(s,0))=(s,\mu \cos(s))$ for $\mu > 0$, we find the following solution:
\begin{align*}
\qquad \varphi=s+\frac{\mu}{2}\sin s(e^{(B_0-1)t}-e^{-(B_0+1)t}),\qquad 
z=B_0t+\frac{\mu}{2}\cos s(e^{(B_0-1)t}+e^{-(B_0+1)t}).
\end{align*}
Note that this solution contains contributions from \eqref{Eq_k_Mode} for all $k$.
It is obvious that this solution does not converge as $t\to\infty$ for any choice of $B_0 \in \R\setminus\{0\}$. 
Again, this is reflected in the critical term in \eqref{energy-equality}
\begin{align*}
\int_0^T\int_{S^1}\Omega(\dot{\gamma},\gamma')dsdt
&= -2\pi B_0^2t + B_0\frac{\pi\mu^2}{4}(e^{-2T(1+B_0)}-e^{-2T(1-B_0)}),
\end{align*}
which is not bounded below with respect to $t$ for any choice of $B_0 \in \R\setminus\{0\}$. Note that the first term on the right hand side is multiplied by $B_0^2$. 
Thus, we cannot change the asymptotic behaviour of the solution by changing the sign of \(B_0\), even if $|B_0|$ is small.
\end{enumerate}

One may of course produce other explicit solutions using the Fourier decomposition of $\xi$ from \eqref{Eq_k_Mode}.

Comparing the findings in this particular example with the general results from Theorem \ref{ThmConvergenceOttarson} and Theorem \ref{ThmConvExactCase}, we can draw the following conclusions:

\begin{Bem} \label{RemTorusEx}
\begin{enumerate}
 \item In general, we cannot  expect to obtain the (sub)convergence of the flow to a magnetic geodesic in case $\Omega$ is not exact unless the initial curve $\gamma_0$ is contractible. In fact, since all 
       solutions in this example are contractible, the initial conditions, being in the same homotopy class, have to be contractible, too. 
       Part (2) of the example shows an initial curve wrapping around the torus and we have seen that the corresponding flow never converges.    
 \item Even if $\Omega$ is exact, we need additional assumptions to guarantee $|A|_{L^\infty} < \infty$ and hence (sub)-convergence of the flow to a magnetic geodesic. In fact, the flow from part 
       \eqref{ExTorusPart2}, viewed as an element of $C^\infty(S^1 \times [0,\infty), C)$, is not subconvergent, which reflects the fact that the potential $A = z d\phi$ is not essentially bounded on the non-compact cylinder $C$. 
 \item An upper bound on $|Z|_{L^\infty}$ (as the one used in Lemma \ref{LemZAssump}) is crucial to obtain convergence. In fact, choosing $k := 1 < B_0$ in part \eqref{ExTorusPart1} above provides an example 
       of a divergent flow. 
 \item Even if the flow is subconvergent, the limit need not be a magnetic geodesic. Setting $\mu := 0$ in part \eqref{ExTorusPart2}, we obtain the constant sequence of curves $\gamma_{t_m} \in 
       C^\infty(S^1,T^2)$ for $t_m := 2\pi m$, $m \in \N$. We have already seen that $\gamma_\infty$ cannot be a magnetic geodesic because it winds one time around the torus.
 \item Provided $|B_0| \leq k$, part \eqref{ExTorusPart1} gives an example, where Proposition \ref{LemNonTriv1} can be applied even though $\Omega$ is not exact. Moreover, we see that the condition 
       $E(\gamma_0) \leq 0$ cannot be improved: Choosing $a=b=1$, we obtain $E(\gamma_0) =  E_{kin}(\gamma_0) + E_A(\gamma_0) = \pi(k^2+kB_0)$. 
       Thus, we find examples of a flow with positive total energy arbitrary close to zero, which shrinks to a point. 
       On the other hand, we also see that the aforementioned condition is not necessary since $B_0 = k$  yields an example with positive initial energy converging to a nontrivial magnetic geodesic. \label{RemTorusExEnergy}
 \item Part \eqref{ExTorusPart1} also illustrates the rescaling technique used in Corollary \ref{CorRescaling}: 
       For parameter $k > 1$, we may consider the flow as map $S^1_{1/k} \times [0,\infty) \rightarrow T^2$. 
       Rescaling it by $\frac{1}{k}$, we get a new flow, defined on $S^1_1 \times [0,\infty)$, which solves \eqref{rescaledEOM} for $\Lambda = \frac{1}{k}$. 
       The rescaled flow converges provided $|\frac{B_0}{k}| \leq 1$ by \ref{RemZ_Eins_Optimal}. 
       This is equivalent to the condition $|B_0| \leq k$, which was obtained above by an explicit calculation. 
       In particular, we find that convergence of the flow to magnetic geodesics may be obtained for arbitrary large values of $|B_0|$ by suitable rescaling in accordance with Corollary \ref{CorRescaling}.
\end{enumerate}
\end{Bem}

\subsection{Magnetic Geodesics on the two-dimensional sphere} \label{ExSphere}
Now we study the case of a two-dimensional spherical target, \(N := S^2\subset\R^3\). 
Magnetic geodesics on \(S^2\) have been extensively studied in \cite{MR2788659}. A natural choice for a magnetic force is given by \(Z(\gamma')=B_0\gamma\times\gamma'\). 
Again, it is not hard to see that the associated 2-form is given by $\Omega = B_0 vol_{S^2}$, which is not exact. In this case, the equation for magnetic geodesics on \(S^2\) acquires the form
\begin{equation}
\label{magnetic-sphere}
\gamma''=-|\gamma'|^2\gamma+B_0\gamma\times\gamma'.
\end{equation}
By a direct computation one can check that
\begin{equation}
\label{solution-sphere}
\gamma(s)=(\sin\theta_0\cos(\frac{B_0}{\cos\theta_0}s),\sin\theta_0\sin(\frac{B_0}{\cos\theta_0}s),\cos\theta_0) 
\end{equation}
solves \eqref{magnetic-sphere} for a fixed value of \(\theta_0\). 
However, note that this equation implicitly assumes that this solution is defined on a circle $S^1_r$ whose radius $r$ is an integer multiple of $\frac{B_0}{\cos\theta_0}$. 
The corresponding heat flow leads to the following system
\begin{equation}
\dot{\gamma}=\gamma''+|\gamma'|^2\gamma-B_0\gamma\times\gamma'
\end{equation}
with initial data \(\gamma_0\).
Here, the relevant term in  \eqref{energy-equality} term can be expressed as
\[
\Omega(\dot{\gamma},\gamma')=B_0\langle\dot{\gamma},\gamma\times\gamma'\rangle=B_0\det(\dot{\gamma},\gamma,\gamma')
\]
and depending on its sign, we expect that the evolution equation converges or not.
Let us analyze this evolution equation for several different ansätze.
\begin{enumerate}
 \item We use the following ansatz:
\begin{equation} \label{EqSphereAnsatz1}
\gamma(s,t) = (\sin\theta(t)\cos\varphi(s),\sin\theta(t)\sin\varphi(s),\cos\theta(t))
\end{equation}
This ansatz leads to a system of differential equations,
which can be simplified easily. From the equation for \(\gamma_3\) 
we obtain
\[
\frac{\dot\theta}{\sin\theta}=B_0\varphi'-\varphi'^2\cos\theta.
\]
Equating the expressions for \(\gamma_1\) and \(\gamma_2\) yields
\(\tan(\varphi)\varphi''=-\cot(\varphi)\varphi''\), which implies that \(\varphi''=0\).
Hence, $\varphi(s) = c\cdot s$ and after a suitable rescaling, which only affects the value of $B_0$, we may assume \(\varphi(s)=s\).  The equations for \(\gamma_i,i=1,2,3\) now lead to
\begin{equation}
\label{evolution-theta-sphere}
\dot{\theta}=-\sin\theta\cos\theta+B_0\sin\theta = \sin\theta (B_0 - \cos\theta). 
\end{equation}
This equation cannot be integrated directly for arbitrary \(B_0\).

To understand the asymptotics of \eqref{evolution-theta-sphere}, we study the zero's of the right hand side.
It can easily be checked that the right hand side vanishes for \(\theta=0,\pi,\theta_0:= \arccos(B_0)\) and the functions $\theta(t)=0,\pi,\theta_0$ are clearly solutions of \eqref{evolution-theta-sphere}. 
Moreover, it is straightforward to check that
\begin{align} \label{EqThetaSigns}
\left.
\begin{array}{c}
\theta > \theta_0 \\
\theta < \theta_0 
\end{array}
\right\}
&\Longrightarrow 
\left\{
\begin{array}{c}
\sin\theta (B_0 - \cos\theta) > 0 \\
\sin\theta (B_0 - \cos\theta) < 0 
\end{array}
\right.
\end{align}
Thus, depending on the initial condition, the solution will converge to \(\theta_\infty=0\), \(\theta_\infty=\pi\) 
or, more interestingly, to \(\theta_\infty=\theta_0\), see also \cite{MR2961944}, Lemma 1.1.
Moreover we see that if $\theta(0)\neq \theta_0$, then the flow equation \eqref{evolution-theta-sphere} causes $\theta(t)$ to approach one of the poles, i.e. $\lim_{t \rightarrow \infty}\theta(t) = 0, \pi$. 
The corresponding limit curve $\gamma_\infty$ will either be one of the poles of the sphere (and hence constant) or given by the curve
\begin{equation}
\gamma_\infty(s)=(\sqrt{1-B_0^2}\cos s,\sqrt{1-B_0^2}\sin s,B_0),
\end{equation}
which is a reparametrization of \eqref{solution-sphere}. We have seen above that unless we already start on this particular solution ($\theta(0) = \theta_0$), the flow will converge to one of the poles. 
Hence, \eqref{solution-sphere} represents a non-stable solution, whereas both trivial solutions are stable, at least within the limits of our ansatz \eqref{EqSphereAnsatz1}.\\
In addition, we check the critical term in \eqref{energy-equality} for this ansatz:
\[
\Omega(\dot{\gamma},\gamma')=B_0\langle\dot{\gamma},\gamma\times\gamma'\rangle=B_0\varphi'\tfrac{d}{dt}(\cos\theta).
\]
When integrating this expression over \(S^1\) we find that it vanishes, in accordance with Remark \ref{RemGeneralConvergence} and the fact that this ansatz converges.

\item We now use a different ansatz:
\begin{equation}
\gamma(s,t)=(\sin\theta(s)\cos\varphi(t),\sin\theta(s)\sin\varphi(t),\cos\theta(s))
\end{equation}
The flow equation now reads \(Z(\gamma')=B_0(-\sin\varphi(t)\theta',-\cos\varphi(t)\theta',0)\) 
and the equation for \(\gamma_3\) then directly implies \(\theta''=0\).
Multiplying the equation for \(\gamma_1\) with \(\sin\varphi\),
multiplying the equation for \(\gamma_2\) with \(\cos\varphi\) and adding up
both contributions then yields
\[
\dot{\varphi}=-B_0\frac{\theta'}{\sin\theta}.
\]
But this means that both sides have to be equal to a constant \(\lambda\), leading to
\[
-\theta'=\lambda\sin\theta.
\]
In addition, we know that \(\theta'\) is constant, such that \(\theta=0\).
Hence, the curve \(\gamma\) will stay at the north pole of \(S^2\).
The critical term in \eqref{energy-equality} takes the form
\[
\Omega(\dot{\gamma},\gamma')=-B_0\dot{\varphi}\tfrac{d}{ds}(\cos\theta)
\]
and this expression again vanishes when integrated over \(S^1\).
\end{enumerate}

In addition to the comments at the end of Example \ref{ExTorus}, we can draw additional conclusions concerning the heat flow approach from the current example on $S^2$: 

\begin{Bem} \label{RemSphere}
\begin{enumerate}
\item By \eqref{solution-sphere}, there exist nontrivial contractible periodic magnetic geodesics. However, the heat flow method is not appropriate to establish their existence since the nontrivial solution is unstable.
Given a generic choice for the initial curve $\gamma_0$, the flow will converge to one of the poles. This behaviour is analogous to the properties of the ordinary heat flow for closed geodesics on $S^2$.
\item The discussion of ansatz (1) also provides an example, where the criterion in Proposition \ref{LemNonTriv1} for non-triviality of $\gamma_\infty$ explicitly fails. 
Choosing a potential $A_p$ on $S^2\setminus \{p\}$ for some point $p \in S^2$ and an initial curve $\gamma_0$ such that $\tfrac{1}{2}\int |\gamma_0'|^2ds + \int \gamma_0^\ast Ads \leq 0$, 
it may be shown that the flow either passes through $p$ at some time or converges to $p$. Thus, the argument in \ref{LemNonTriv1} breaks down.   
\end{enumerate} 
\end{Bem}

\subsection{Magnetic Geodesics on the hyperbolic plane} \label{ExHyperbolic}
Finally, we assume that \(N\) is the hyperbolic plane \(\hyp^2\) with constant curvature \(-1\).
A general existence result for magnetic geodesics on \(\hyp^2\) was recently established in \cite{MR2959932}.
We choose the convention
\[
\hyp^2 = \{ (\gamma_1,\gamma_2,\gamma_3) \in \R^3 \mid \gamma_1^2-\gamma_2^2-\gamma_3^2=1 \},
\]
i.e the axis of rotation of the hyperboloid points in the $\gamma_1$-direction. 
As in \cite{MR2959932}, p.6, we define a modified cross product between two vectors \(v,w\) by
\begin{equation}
v\tilde{\times}w=(v_3w_2-v_2w_3,v_1w_3-v_3w_1,v_1w_2-v_2w_1).
\end{equation}
The equation for magnetic geodesics is now given by
\[
\gamma''=\gamma|\gamma'|^2+B_0\gamma\tilde{\times}\gamma'
\]
and as in the spherical case, an explicit solution can easily be found:
\begin{equation}
\label{solution-hyperbolic}
\gamma(s)=(\cosh\theta_0,\sinh\theta_0\cos(\frac{B_0}{\cosh\theta_0}s),\sinh\theta_0\sin(\frac{B_0}{\cosh\theta_0}s)). 
\end{equation}
Note that for \(B_0\to 0\), the solution shrinks to a point. This reflects the fact that there do not exist ordinary closed geodesics on \(\hyp^2\).\\

The evolution equation for magnetic geodesics then reads
\begin{equation}
\dot{\gamma}=\gamma''-\gamma|\gamma'|^2-B_0\gamma
\end{equation}
for some given initial data \(\gamma_0\). Let us make the following ansatz
\[
\gamma(s,t)=(\cosh\theta(t),\sinh\theta(t)\cos\varphi(s),\sinh\theta(t)\sin\varphi(s)).
\]
Similar to the spherical case this leads to the following two equations
\begin{align*}
\dot\theta=-\cosh\theta\sinh\theta\varphi'^2-\sinh\theta\varphi',\qquad \varphi''=0.
\end{align*}
Thus, setting \(\varphi=s\), we obtain
\[
\dot\theta=-\cosh\theta\sinh\theta+B_0\sinh\theta.
\]
As in the spherical case, we cannot integrate this equation directly. 
Thus, we again analyze the zero's of the right hand side, which are given by \(0,\theta_0 := \operatorname{arcosh} B_0\).
Similar to \eqref{EqThetaSigns}, we moreover have
\begin{align} \label{EqThetaSignsHyp}
\left.
\begin{array}{c}
\theta > \theta_0 \\
\theta < \theta_0 
\end{array}
\right\}
&\Longrightarrow 
\left\{
\begin{array}{c}
\sinh\theta (B_0 - \cosh\theta) < 0 \\
\sinh\theta (B_0 - \cosh\theta) > 0 
\end{array}
\right.
\end{align}
In contrast to the spherical case, we now find that unless $\theta(0) = 0$, we have $\lim_{t\rightarrow \infty}\theta(t) = \theta_0$.
Hence, depending on the initial condition, the solution will stay at $(1,0,0)$ if it started at this point or converges to the curve
\begin{equation}
\gamma_\infty=(B_0,\sqrt{B_0^2-1}\cos s,\sqrt{B_0^2-1}\sin s),
\end{equation}
which again coincides with \eqref{solution-hyperbolic} up to a reparametrization. 

\begin{Bem}
In contrast to Example \ref{ExSphere}, the nontrivial solution on $\hyp^2$ is stable, at least within the limits of our ansatz. 
In particular, the heat flow method is appropriate for establishing the existence of the nontrivial magnetic geodesics \eqref{solution-hyperbolic}. 
Again, this reflects the properties of the heat flow for ordinary geodesics where critical points are stable if the manifold has strictly negative curvature.
Due to the magnetic terms in \eqref{second-variation}, we do not fully understand the role of curvature for magnetic geodesics at present.     
\end{Bem}

\bibliographystyle{alpha}
\bibliography{mybib}
\end{document}